\documentclass[10pt]{article}
\usepackage{amssymb}
\usepackage{amsmath,amssymb}\usepackage{cite}\usepackage{mathrsfs}\usepackage[amsmath,thmmarks]{ntheorem}
\usepackage{listings}
\usepackage[titletoc]{appendix}\allowdisplaybreaks

\makeatletter
\renewcommand\theequation{\thesection.\@arabic\c@equation}
\makeatother

\newtheorem{thm}{ Theorem}[section]%
\newtheorem{lem}[thm]{ Lemma}%
\newtheorem{cor}[thm]{ Corollary}%

\newtheorem{Remark}[thm]{ Remark}%
\newtheorem{Pro}[thm]{ Proposition}%
\newtheorem{Que}[thm]{Question}%
\newtheorem{Fac}[thm]{Fact}%
\input cyracc.def

\def\f{\noindent}

\def\demo{\f{\bf Proof}\hskip10pt}

\def\qed{\hfill $\Box$}

\begin{document}
\title{\textbf{The full automorphism groups of general position graphs\\{{\small Dedicated to my father Hongqi Pan's 75th birthday}}\thanks{The research of the work was partially supported
by the Hainan Provincial Natural Science Foundation of China (No. 122RC652) and the National Natural Science Foundation of China (No. 12061030; No. 61962018).}}}

\footnotetext{Junyao Pan. E-mail addresses: $\rm{Junyao_{-}Pan}$@126.com}

\author{Junyao Pan\,
 \\\\
School of Sciences, University of Wuxi, Wuxi, Jiangsu,\\ 214105
 People's Republic of China \\}
\date {} \maketitle

\baselineskip=16pt

\vskip0.5cm

{\bf Abstract:} Let $S$ be a non-empty finite set. A flag of $S$ is a set $f$ of non-empty proper subsets of $S$ such that $X\subseteq Y$ or $Y\subseteq X$ for all $X,Y\in f$. The set $\{|X|:X\in f\}$ is called the type of $f$. Two flags $f$ and $f'$ are in general position with respect to $S$ if $X\cap Y=\emptyset$ or $X\cup Y=S$ for all $X\in f$ and $Y\in f'$. For a fixed type $T$, Klaus Metsch defined the general position graph $\Gamma(S,T)$ whose vertices are the flags of $S$ of type $T$ with two vertices being adjacent when the corresponding flags are in general position. In this paper, we characterize the full automorphism groups of $\Gamma(S,T)$ in the case that $|T|=2$. In particular, we solve an open problem proposed by Klaus Metsch.



\section {Introduction}
Throughout this paper, $[n]=\{1,2,...,n\}$ denotes the standard $n$-element set and $\overline{A}$ stands for the complement of $A$ in $[n]$ where $A\subseteq[n]$. Moreover, $C^k_m=\frac{m(m-1)\cdot\cdot\cdot(m-k+1)}{k!}$ expresses the binomial coefficient.

For two positive integers $n$ and $k$ with $n\geq2k$, the \emph{Kneser graph} $KG(n,k)$ has as vertices the $k$-subsets of $[n]$ with edges defined by disjoint pairs of $k$-subsets. It is well-known that the problem of the independence number of $KG(n,k)$ reduces to the famous Erd\rm{\H{o}}s-Ko-Rado Theorem. From this perspective, Klaus Metsch \cite{M} introduced the \emph{general position graph} of flags of $[n]$ to generalize the famous Erd\rm{\H{o}}s-Ko-Rado Theorem. A \emph{flag} of $[n]$ is a set $f$ of non-empty proper subsets of $[n]$ such that $A\subseteq B$ or $B\subseteq A$ for all $A,B\in f$. The set $\{|A|:A\in f\}$ is called the type of $f$. Two flags $f$ and $f'$ are in \emph{general position} with respect to $[n]$ if $A\cap B=\emptyset$ or $A\cup B=[n]$ for all $A\in f$ and $B\in f'$. For a fixed type $T\subseteq[n-1]$ the \emph{general position graph} whose vertices are the flags of $[n]$ of type $T$ with two vertices being adjacent when the corresponding flags are in general position will be abbreviated as $\Gamma(n,T)$. If $|T|=1$, then this graph is isomorphic to a corresponding Kneser graph. Klaus Metsch \cite{M} not only described the independence number of $\Gamma(n,T)$ in some situations but also proposed several interesting open problems, such as the following one:

\begin{Que}\label{pan1-1}\normalfont(\cite[Problem~5]{M}~)
Is it true that the graphs $\Gamma(n,T)$ have $S_n$ as automorphism group, where $T=\{a,b\}$ with $n\geq a+b+1$ and $a<\frac{n}{2}<b$? Can \cite[Remark~5.13]{M} be used to show this?
\end{Que}

Let $\Gamma = (V, E)$ be a undirected graph with vertex set $V$ and edge set $E$. If there exists a bijection $\alpha$ from $V$ to $V$ such that $(f^\alpha,g^\alpha)\in E$ if and only if $(f,g)\in E$ for all $f,g\in V$, then $\alpha$ is called an automorphism of $\Gamma$. Let $Aut(\Gamma)$ denote the full automorphism group of $\Gamma$. Actually, the research on the automorphism groups of graphs has always been an interesting topic for many scholars in group theory and graph theory, for examples \cite{F,FMW,G,KNZ,PG,T}. Thereby, we are interested in Question\ \ref{pan1-1}.

Review some notions and notations about permutation groups, for details see \cite{C,D}. Let $G$ be a transitive permutation group acting on $[n]$. A non-empty subset $\Delta$ of $[n]$ is called a \emph{block} for $G$ if for each $\alpha\in G$ either $\Delta^\alpha=\Delta$ or $\Delta^\alpha\cap\Delta=\emptyset$. Clearly, the singletons $\{i\}$ $(i\in[n])$ and $[n]$ are blocks, and so these blocks are called the \emph{trivial} blocks. Any other block is called \emph{nontrivial}. Put $\Sigma=\{\Delta^\alpha:\alpha\in G\}$ where $\Delta$ is a block of $G$. We call $\Sigma$ the system of blocks containing $\Delta$. Clearly, $G$ reduces a permutation group acting on $\Sigma$, denoted by $G|_\Sigma$. In addition, there exists a natural homomorphism from $G$ to $G|_\Sigma$, and the kernel of this homomorphism consists of all permutations in $G$ which fix every block in $\Sigma$. In this note, we divide three cases to construct systems of blocks of $Aut(\Gamma(n,T))$ acting flags and further we show the kernels of the corresponding homomorphisms are all trivial. Thus, we give a positive answer to the Question\ \ref{pan1-1}.

\section {Preliminaries}

Let $\mathcal{F}^T_n$ denote the set of all flags of type $T$ of $[n]$. In other words, $\mathcal{F}^T_n$ is the vertex set of $\Gamma(n,T)$. Moreover, we set $\mathcal{F}^{(T|A)}_n=\{f\in\mathcal{F}^T_n:~A\in f\}$, where $A\subseteq[n]$ and $|A|\in T$. Let $f$ and $g$ be two flags in $\mathcal{F}^T_n$. If there exists an edge between $f$ and $g$ in $\Gamma(n,T)$, then $f$ and $g$ are called \emph{neighbor} (see \cite{R}). In addition, $N(f)$ stands for the collection of all neighbours of vertex $f$. Here, we state a well-known fact that is the key idea of solving Question\ \ref{pan1-1}.

\begin{Fac}\label{pan2-0}\normalfont
Let $f,g\in\mathcal{F}^T_n$. Then $|N(f)^\alpha|=|N(f^\alpha)|$ and $|N(f)\cap N(g)|=|N(f^\alpha)\cap N(g^\alpha)|$ for every $\alpha\in Aut(\Gamma(n,T))$.
\end{Fac}
For convenience, we set $N(f,g)=N(f)\cap N(g)$ and $N_{m}(n|T)={\rm{max}}\{|N(f,g)|:f,g\in\mathcal{F}^T_n\}$. Next, we characterize $N_m(n|T)$ in some situations.

\begin{Pro}\label{pan2-1}\normalfont
Let $f=\{A,B\}$ and $g=\{C,D\}$ be two flags in $\mathcal{F}^T_n$ with $|A|=|C|=a$ and $|B|=|D|=b$, where $T=\{a,b\}$ with $a+b+1\leq n$ and $a<\frac{n}{2}<b$. Then the followings hold.

(i) If $b>\frac{2n}{3}$, then $|N(f,g)|=N_m(n|T)$ if and only if $A=C$ and $|B\cap D|=b-1$. In this case, $N_m(n|T)=C^a_{n-b-1}C^{2b-n-1}_{b-a-1}$.

(ii) If $b<\frac{2n}{3}$, then $|N(f,g)|=N_m(n|T)$ if and only if $B=D$ and $|A\cap C|=a-1$. In this case, $N_m(n|T)=C^a_{n-b}C^{2b-n}_{b-a-1}$.

(iii) If $b=\frac{2n}{3}$, then $|N(f,g)|=N_m(n|T)$ if and only if either $A=C$, $|B\cap D|=b-1$ or $B=D$, $|A\cap C|=a-1$. In this case, $N_m(n|T)=C^a_{n-b-1}C^{2b-n-1}_{b-a-1}=C^a_{n-b}C^{2b-n}_{b-a-1}$.
\end{Pro}

\demo Let $h=\{G,H\}\in\mathcal{F}^T_n$ with $|G|=a$ and $|H|=b$. Since $a+b+1\leq n$ and $a<\frac{n}{2}<b$, it follows that $h\in N(f,g)$ if and only if $H\cap A=\emptyset$, $H\cap C=\emptyset$, $H\cup B=H\cup D=[n]$, $G\cap B=\emptyset$ and $G\cap D=\emptyset$. In other words, $h\in N(f,g)$ if and only if $G\subseteq[n]\setminus(B\cup D)$ and $\overline{B\cap D}\subseteq H$ and $H\cap (A\cup C)=\emptyset$. Clearly, the number of the choices of $G$ is $C^a_{n-|B\cup D|}$. Consider $H$. Due to $\overline{B\cap D}\subseteq H$ and $H\cap(A\cup C)=\emptyset$, it follows that $A\cup C\subseteq B\cap D$ and otherwise $N(f,g)=\emptyset$. Note that the number of the choices of $H$ is equal to the number of the choices of $(b-|\overline{B\cap D}|)$-subsets from $(B\cap D)\setminus(A\cup C)$. Therefore, the number of the choices of $H$ is $C^{b-|\overline{B\cap D}|}_{|B\cap D|-|A\cup C|}$. Additionally, it is clear that $|B\cup D|=2b-|B\cap D|$ and $|\overline{B\cap D}|=n-|B\cap D|$. Thus, we have deduced that
\begin{equation}\label{eq111}
|N(f,g)|=C^a_{n-2b+|B\cap D|}C^{b+|B\cap D|-n}_{|B\cap D|-|A\cup C|}=C^a_{n-2b+|B\cap D|}C^{n-b-|A\cup C|}_{|B\cap D|-|A\cup C|}.
\end{equation}
Here, we make a statement that every abnormal situation occurs imply $N(f,g)=\emptyset$, where the abnormal situation is $|B\cap D|-|A\cup C|<0$ or $n-2b+|B\cap D|<a$ and so on. Our goal is to find $N_m(n|T)$, and so we do not discuss when $N(f,g)=\emptyset$ holds.

Fixing $|A\cup C|\geq a$. The equality\ \ref{eq111} shows that $|N(f,g)|$ increases with the increase of $|B\cap D|$. Likewise, fixing $|B\cap D|\leq b$, $|N(f,g)|$ decreases with the increase of $|A\cup C|$. Therefore, $N_m(n|T)$ occurs in two possible situations. One is that $|A\cup C|=a$ and $|B\cap D|=b-1$, and the other is that $|A\cup C|=a+1$ and $|B\cap D|=b$. If $|A\cup C|=a$ and $|B\cap D|=b-1$, then by equality\ \ref{eq111} we deduce that
\begin{equation}\label{eq11}
|N(f,g)|=C^a_{n-b-1}C^{2b-n-1}_{b-a-1}.
\end{equation}
If $|A\cup C|=a+1$ and $|B\cap D|=b$, then by equality\ \ref{eq111} we infer that
\begin{equation}\label{eq2}
|N(f,g)|=C^a_{n-b}C^{2b-n}_{b-a-1}.
\end{equation}
Compare the equalities (\ \ref{eq11}~) and (\ \ref{eq2}~). We see that $\frac{(\ \ref{eq11}~)}{(\ \ref{eq2}~)}=\frac{C^a_{n-b-1}C^{2b-n-1}_{b-a-1}}{C^a_{n-b}C^{2b-n}_{b-a-1}}=\frac{2b-n}{n-b}$. Obviously, if $b>\frac{2n}{3}$, then $\frac{\ \ref{eq11}}{\ \ref{eq2}}>1$ and so $N_m(n|T)=C^a_{n-b-1}C^{2b-n-1}_{b-a-1}$. In this case, $|N(f,g)|=N_m(n|T)$ if and only if $A=C$ and $|B\cap D|=b-1$. If $b<\frac{2n}{3}$, then $\frac{\ \ref{eq11}}{\ \ref{eq2}}<1$ and thus $N_m(n|T)=C^a_{n-b}C^{2b-n}_{b-a-1}$. In this case, $|N(f,g)|=N_m(n|T)$ if and only if $B=D$ and $|A\cap C|=a-1$. Additionally, if $b=\frac{2n}{3}$ then $\frac{\ \ref{eq11}}{\ \ref{eq2}}=1$ and thus $|N(f,g)|=N_m(n|T)$ if and only if either $A=C$ and $|B\cap D|=b-1$ or $B=D$ and $|A\cap C|=a-1$. In this case, $N_m(n|T)=C^a_{n-b-1}C^{2b-n-1}_{b-a-1}=C^a_{n-b}C^{2b-n}_{b-a-1}$.     \qed

\begin{Remark}\label{pan2-2}\normalfont
Let $f,g\in\mathcal{F}^T_n$ such that $|N(f,g)|=N_{m}(n|T)$. Then $|N(f^\alpha,g^\alpha)|=N_{m}(n|T)$ for any $\alpha\in Aut(\Gamma(n,T))$.
\end{Remark}
\demo This remark follows from Fact\ \ref{pan2-0}.    \qed

By Proposition\ \ref{pan2-1} (iii), we see that if $b=\frac{2n}{3}$ then $N_{m}(n|T)$ occurs in two cases. This urges us to further study this situation. Now we investigate the second maximum of $|N(f,g)|$ and so we define $N_{sm}(n|T)={\rm{max}}\{|N(f,g)|:f,g\in\mathcal{F}^T_n,|N(f,g)|<N_{m}(n|T)\}$.

\begin{Pro}\label{pan2-3}\normalfont
Let $f=\{A,B\}$ and $g=\{C,D\}$ be two flags in $\mathcal{F}^T_n$ with $|A|=|C|=a$ and $|B|=|D|=b$, where $T=\{a,b\}$ with $a\leq\frac{n}{3}-1$ and $b=\frac{2n}{3}$. Then, $|N(f,g)|=N_{sm}(n|T)$ if and only if $|A\cup C|=a+1$ and $|B\cup D|=b+1$.
\end{Pro}

\demo According to the equality \ \ref{eq111}, we see that $|N(f,g)|=N_{sm}(n|T)$ occurs in three possible cases, those are, $|A\cup C|=a+1$, $|B\cup D|=b+1$ or $|A\cup C|=a+2$, $|B\cup D|=b$ or $|A\cup C|=a$, $|B\cup D|=b+2$. Next, we start to compute $|N(f,g)|$ in three cases respectively.

If $|A\cup C|=a+1$ and $|B\cup D|=b+1$, then $|B\cap D|=b-1$. In this case, by equality\ \ref{eq111}, we infer that
\begin{equation}\label{eq3}
|N(f,g)|=C^a_{n-b-1}C^{2b-n-1}_{b-a-2}.
\end{equation}
Similarly, if $|A\cup C|=a$ and $|B\cup D|=b+2$, then
\begin{equation}\label{eq4}
|N(f,g)|=C^a_{n-b-2}C^{2b-n-2}_{b-a-2};
\end{equation}
and if $|A\cup C|=a+2$ and $|B\cup D|=b$ then
\begin{equation}\label{eq5}
|N(f,g)|=C^a_{n-b}C^{2b-n}_{b-a-2}.
\end{equation}
Compare three equalities (\ \ref{eq3}~) and (\ \ref{eq4}~) and (\ \ref{eq5}~). We deduce that $$\frac{(\ \ref{eq4}~)}{(\ \ref{eq3}~)}=\frac{C^a_{n-b-2}C^{2b-n-2}_{b-a-2}}{C^a_{n-b-1}C^{2b-n-1}_{b-a-2}}=\frac{(n-b-a-1)(2b-n-1)}{(n-b-1)(n-b-a)}=\frac{n-b-a-1}{n-b-a}<1~\rm{and}$$
$$\frac{(\ \ref{eq5}~)}{(\ \ref{eq3}~)}=\frac{C^a_{n-b}C^{2b-n}_{b-a-2}}{C^a_{n-b-1}C^{2b-n-1}_{b-a-2}}=\frac{(n-b)(n-b-a-1)}{(n-b-a)(2b-n)}=\frac{n-b-a-1}{n-b-a}<1.$$
Therefore, $|N(f,g)|=N_{sm}(n|T)$ if and only if $|A\cup C|=a+1$ and $|B\cup D|=b+1$. In particular, $N_{sm}(n|T)=C^a_{n-b-1}C^{2b-n-1}_{b-a-2}$.    \qed

So far, we have seen that $N_{sm}(n|T)$ occurs in unique form when $b=\frac{2n}{3}$.
Let $f,g\in\mathcal{F}^T_n$ where $T=\{a,b\}$ with $a\leq\frac{n}{3}-1$ and $b=\frac{2n}{3}$. Define $SM(f)=\{h\in\mathcal{F}^T_n:|N(f,h)|=N_{sm}(n|T)\}$ and $SM(f,g)=SM(f)\cap SM(g)$. Next we state two results which are useful in dealing with the case that $b=\frac{2n}{3}$.

\begin{Pro}\label{pan2-4}\normalfont
Suppose that $f,g\in\mathcal{F}^T_n$, where $T=\{a,b\}$ with $a\leq\frac{n}{3}-1$ and $b=\frac{2n}{3}$. Then $|SM(f,g)|=|SM(f^\alpha,g^\alpha)|$ for any $\alpha\in Aut(\Gamma(n,T))$.
\end{Pro}
\demo Assume that $h\in\mathcal{F}^T_n$ such that $h\in SM(f,g)$. Thus, $|N(f,h)|=|N(g,h)|=N_{sm}(n|T)$. Then by Fact\ \ref{pan2-0} we deduce that $|N(f,h)|=|N(f^\alpha,h^\alpha)|$ and $|N(g,h)|=|N(g^\alpha,h^\alpha)|$ for any $\alpha\in Aut(\Gamma(n,T))$. Conversely, for any $\alpha\in Aut(\Gamma(n,T))$, if $v\in SM(f^\alpha,g^\alpha)$ then $v^{\alpha^{-1}}\in SM(f,g)$. This completes the proof of this proposition.   \qed

\begin{Pro}\label{pan2-5}\normalfont
Let $f=\{A,B\}$, $g=\{A,C\}$, $x=\{D,E\}$ and $y=\{F,E\}$ be four flags in $\mathcal{F}^T_n$ such that $|B\cap C|=b-1$ and $|D\cap F|=a-1$, where $T=\{a,b\}$ with $a\leq\frac{n}{3}-1$ and $b=\frac{2n}{3}$. Then $|SM(f,g)|\neq|SM(x,y)|$.
\end{Pro}
\demo Let $h=\{G,H\}\in\mathcal{F}^T_n$. By applying Proposition\ \ref{pan2-3}, we deduce that $h\in SM(f,g)$ if and only if $|G\cap A|=a-1$ and $|H\cap B|=b-1$ and $|H\cap C|=b-1$. Count $|SM(f,g)|$. Note that there exist two possible shapes for $H$, those are, $B\cap C\subseteq H$ and $B\cap C\not\subseteq H$ respectively. In the case of $B\cap C\subseteq H$, we see that $H=(B\cap C)\cup\{i\}$ where $i\in[n]\setminus(B\cup C)$ and so the number of the choices of $H$ is $C^1_{n-b-1}$. Fix an $H$, if $G\subseteq B\cap C$ then $G$ is the union of a $(a-1)$-subset in $A$ and a $1$-subset in $(B\cap C)\setminus A$, and thus the number of the choices of $G$ is $C^{a-1}_{a}C^1_{b-1-a}$; and if $G\not\subseteq B\cap C$ then $G$ is the union of an $(a-1)$-subset in $A$ and $H\setminus(B\cap C)$ and so the number of the choices of $G$ is $C^{a-1}_{a}$. Hence, in the case of $B\cap C\subseteq H$, the number of the choices of $h$ in $SM(f,g)$ is
$$C^1_{n-b-1}(C^{a-1}_{a}C^1_{b-1-a}+C^{a-1}_{a})=(\frac{n}{3}-1)a(\frac{2n}{3}-a-1)+(\frac{n}{3}-1)a.$$
Consider $B\cap C\not\subseteq H$. In this case, $H=(B\setminus C)\cup(C\setminus B)\cup K$ where $K$ is a $(b-2)$-subset of $B\cap C$. In addition, there exist two possible subcases for $H$, those are, $A\subseteq H$ and $A\not\subseteq H$. In the subcase of $A\subseteq H$, the number of the choices of $H$ is $C^{b-2-a}_{b-1-a}$. Fix an $H$, $G$ is the union of an $(a-1)$-subset of $A$ and a $1$-subset of $H\setminus A$, and thus the number of the choices of $G$ is $C^{a-1}_{a}C^1_{b-a}$. If $A\not\subseteq H$, then $H=(B\cup C)\setminus\{j\}$ where $j\in A$, and so the number of the choices of $H$ is $C^{1}_{a}$. Fix an $H$, $G=(A\setminus\{j\})\cup\{k\}$ where $k\in H\setminus A$ and so the number of the choices of $G$ is $C^1_{b-a+1}$. Therefore, in the case of $B\cap C\not\subseteq H$, the number of the choices of $h$ in $SM(f,g)$ is $$C^{b-2-a}_{b-1-a}C^{a-1}_{a}C^1_{b-a}+C^{1}_{a}C^1_{b-a+1}=(\frac{2n}{3}-a-1)a(\frac{2n}{3}-a)+a(\frac{2n}{3}-a+1).$$ So we deduce that
$$|SM(f,g)|=a(\frac{2n}{3}-a)(n-a-1)+a.$$

Likewise, we can count $|SM(x,y)|$. Let $z=\{U,V\}\in\mathcal{F}^T_n$. By applying Proposition\ \ref{pan2-3}, we deduce that $z\in SM(x,y)$ if and only if $|U\cap D|=a-1$ and $|U\cap F|=a-1$ and $|V\cap E|=b-1$. Note that there exist two possible shapes for $U$, those are, $U\subseteq E$ and $U\not\subseteq E$ respectively. Consider $U\subseteq E$. If $D\cap F\subseteq U$, then $U=(D\cap F)\cup\{p\}$ where $p\in E\setminus(D\cup F)$, and so the number of the choices of $U$ is $C^{1}_{b-1-a}$. Fix an $U$, $V$ is the union of a $(b-1)$-subset containing $U$ of $E$ and a $1$-subset of $[n]\setminus E$, and thus the number of the choices of $V$ is $C^{b-a-1}_{b-a}C^1_{n-b}$. If $D\cap F\not\subseteq U$, then $U=(D\cup F)\setminus\{q\}$ where $q\in D\cap F$ and so the number of the choices of $U$ is $C^{1}_{a-1}$. Fix an $U$, the number of the choices of $V$ is $C^{b-a-1}_{b-a}C^1_{n-b}$ too. Therefore, in the case of $U\subseteq E$, the number of the choices of $z$ in $SM(x,y)$ is
$$C^{1}_{b-1-a}C^{b-a-1}_{b-a}C^1_{n-b}+C^{1}_{a-1}C^{b-a-1}_{b-a}C^1_{n-b}=(\frac{2n}{3}-a-1)(\frac{2n}{3}-a)\frac{n}{3}+(a-1)(\frac{2n}{3}-a)\frac{n}{3}.$$
Consider $U\not\subseteq E$. In this case, $U=(D\cap F)\cup\{r\}$ where $r\in[n]\setminus E$ and so the number of the choices of $U$ is $C^{1}_{n-b}$. Fix an $U$, $V$ is the union of $U$ and a $(b-a)$-subset of $E\setminus(D\cap F)$ and thus the number of the choices of $V$ is $C^{b-a}_{b-a+1}$. Hence, in the case of $U\not\subseteq E$, the number of the choices of $z$ in $SM(x,y)$ is $$C^{1}_{n-b}C^{b-a}_{b-a+1}=\frac{n}{3}(\frac{2n}{3}-a+1).$$
Therefore, we deduce that
$$|SM(x,y)|=\frac{n}{3}(\frac{2n}{3}-a)(\frac{2n}{3}-1)+\frac{n}{3}.$$
By computing, we derive $|SM(f,g)|-|SM(x,y)|=(a-\frac{n}{3})\{1-(\frac{2n}{3}-a)(a-\frac{2n}{3}+1)\}\neq0$ for all $a\leq\frac{n}{3}-1$. The proof of this proposition is complete.      \qed

Now we end this section by introducing a graph that will be used. Let $\dbinom{[n]}{k}$ denote the set of all $k$-subsets of $[n]$, where $1\leq k\leq n-1$. If $A,B\in\dbinom{[n]}{k}$ such that $|A\cap B|=k-1$, then we say that $A$ and $B$ are \emph{almost identical}. Define \emph{almost identical graph} $AIG(n,k)$ has as vertex set $\dbinom{[n]}{k}$ with two vertices being adjacent when the corresponding subsets are almost identical.

\begin{Pro}\label{pan2-6}\normalfont
Let $n,k$ be positive integers with $1\leq k\leq n-1$. Then $AIG(n,k)$ is a connected graph.
\end{Pro}
\demo Let $A$ and $B$ be two vertices of $AIG(n,k)$. If $|A\cap B|=k-1$ then $A$ and $B$ are adjacent. So we assume that $A=\{a_1,a_2,...,a_l,b_1,b_2,...,b_{k-l}\}$ and $B=\{a_1,a_2,...,a_l,c_1,c_2,...,c_{k-l}\}$ with $|A\cap B|=|\{a_1,a_2,...,a_l\}|=l<k-1$. We note the characteristic that
$$\{b_1,b_2,...,b_{k-l}\}\leftrightarrow \{b_1,b_2,...,b_{k-l-1},c_1\}\leftrightarrow\{b_1,b_2,...,b_{k-l-2},c_1,c_2\}\leftrightarrow\cdot\cdot\cdot\leftrightarrow\{b_1,c_1,c_2,...,c_{k-l-1}\},$$
which implies that there exists a path between $A$ and $B$. The proof of this proposition is complete.  \qed

\section {Main Result}

It is well-known that if $|T|=1$ then $\Gamma(n,T)$ is isomorphic to a Kneser graph. However, we note that if $T=\{t\}$ with $t>\frac{n}{2}$ then $\Gamma(n,T)$ does not fit the traditional definition of Kneser graph. In order to facilitate the readers, we first give a proof of this well-known result.

\begin{lem}\label{pan3-1}\normalfont
Let $T=\{t\}\subseteq[n-1]$. Then $Aut(\Gamma(n,T))\cong S_n$.
\end{lem}

\demo Let $T=\{t\}\subseteq[n-1]$. If $t\leq\frac{n}{2}$, then by \cite[Corollary~7.8.2]{GR} we see that $ Aut(\Gamma(n,T))\cong S_n$. Suppose $t>\frac{n}{2}$. Define a map $\sigma:\Gamma(n,T)\rightarrow KG(n,n-t), A\mapsto\overline{A}$, where $\overline{A}$ denotes the complement of $A$ in $[n]$. It is easy to verify that $\sigma$ is an isomorphic map from $\Gamma(n,T)$ to $KG(n,n-t)$. The proof of this lemma is complete. \qed

\begin{Pro}\label{pan3-2}\normalfont
Let $T=\{a,b\}\subseteq[n-1]$ such that $a<\frac{n}{2}$ and $\frac{2n}{3}<b$ and $a+b+1\leq n$. Then $\Sigma=\{\mathcal{F}^{(T|A)}_n: A\subseteq[n],|A|=a\}$ is a system of nontrivial blocks of $Aut(\Gamma(n,T))$ acting on $\mathcal{F}^{T}_n$.
\end{Pro}
\demo Let $\rho$ be an automorphism in $Aut(\Gamma(n,T))$. Pick $f=\{A,B\}$ and $g=\{A,C\}$ in $\mathcal{F}^{(T|A)}$ such that $|B\cap C|=b-1$ and $|A|=a$. It follows from Proposition\ \ref{pan2-1} (i) and Remark\ \ref{pan2-2} that there exists an $a$-subset $D\subseteq[n]$ such that $f^\rho,g^\rho\in\mathcal{F}^{(T|D)}$. In particular, if $f^\rho=\{D,G\}$ and $g^\rho=\{D,H\}$ then $|G\cap H|=b-1$. By the same token, we deduce that $h^\rho\in\mathcal{F}^{(T|D)}$ in case when $h=\{A,X\}\in\mathcal{F}^{(T|A)}$ with $|X\cap B|=b-1$ or $|X\cap C|=b-1$. Continue moving forward along this line of thought, Proposition\ \ref{pan2-6} indicates that $h^\rho\in\mathcal{F}^{(T|D)}$ for every $h\in\mathcal{F}^{(T|A)}$. Therefore, $\sum$ is a system of nontrivial blocks of $Aut(\Gamma(n,T))$ acting on $\mathcal{F}^{T}_n$.  \qed

\begin{lem}\label{pan3-3}\normalfont
Let $T=\{a,b\}\subseteq[n-1]$ such that $a<\frac{n}{2}$ and $\frac{2n}{3}<b$ and $a+b+1\leq n$. Then $Aut(\Gamma(n,T))\cong S_n$.
\end{lem}
\demo It follows from Proposition\ \ref{pan3-2} that $\sum=\{\mathcal{F}^{(T|A)}_n: A\subseteq[n],|A|=a\}$ is a system of blocks of $Aut(\Gamma(n,T))$ acting on $\mathcal{F}^{T}_n$. Consider the induced action of $Aut(\Gamma(n,T))$ acting on $\sum$, that is, $Aut(\Gamma(n,T))|_{\sum}$. Let $\rho\in Aut(\Gamma(n,T))$ and ${\mathcal{F}}^{(T|A)}_n,{\mathcal{F}}^{(T|B)}_n,{\mathcal{F}}^{(T|C)}_n,{\mathcal{F}}^{(T|D)}_n\in\Sigma$ such that ${{\mathcal{F}}^{(T|A)}_n}^{\rho}={{\mathcal{F}}^{(T|B)}_n}$ and ${{\mathcal{F}}^{(T|C)}_n}^{\rho}={{\mathcal{F}}^{(T|D)}_n}$. Note that $B\cap D=\emptyset$ if and only if $A\cap C=\emptyset$. Otherwise one of $\{(f,g),(f^\rho,g^\rho)\}$ is an edge and the other is not an edge for some $f\in\mathcal{F}^{(T|A)}_n$ and $g\in\mathcal{F}^{(T|C)}_n$. Therefore, $Aut(\Gamma(n,T))|_{\Sigma}\leq Aut(KG(n,a))$. On the other hand, the symmetric group on the set $[n]$ induces an automorphism group of $\Gamma(n,T)$, and thus $Aut(\Gamma(n,T))|_{\Sigma}\cong Aut(KG(n,a))\cong S_n$. Thus, it suffices to check that no nonidentity automorphism can fix all the blocks of $\sum$. Let $\rho$ be an automorphism in $Aut(\Gamma(n,T))$ such that ${{\mathcal{F}}^{(T|A)}_n}^{\rho}={{\mathcal{F}}^{(T|A)}_n}$ for all $A\subseteq[n]$ with $|A|=a$. Assume that there exist two distinct flags $f=\{A,B\}$ and $g=\{A,C\}$ in $\mathcal{F}^{(T|A)}_n$ such that $f^\rho=g$ for some $A\subseteq[n]$ with $|A|=a$. Since $f\neq g$, it follows that $B\neq C$, in other words, there exists an $i\in C\setminus B$. Hence, there exists a flag $h=\{D,E\}\in N(f)$ such that $i\in D$ and $|D|=a$. By Fact\ \ref{pan2-0}, we see that $N(f)^\rho=N(g)$. However, it is clear that $N(g)\cap\mathcal{F}^{(T|D)}_n=\emptyset$, and which indicates that ${{\mathcal{F}}^{(T|D)}_n}^{\rho}\neq{{\mathcal{F}}^{(T|D)}_n}$, a contradiction.    \qed

\begin{Pro}\label{pan3-4}\normalfont
Let $T=\{a,b\}\subseteq[n-1]$ such that $a<\frac{n}{2}<b<\frac{2n}{3}$ and $a+b+1\leq n$. Then $\Omega=\{\mathcal{F}^{(T|B)}_n: B\subseteq[n],|B|=b\}$ is a system of nontrivial blocks of $Aut(\Gamma(n,T))$ acting on $\mathcal{F}^{T}_n$.
\end{Pro}
\demo Let $\rho$ be an automorphism in $Aut(\Gamma(n,T))$. Take two flags $f=\{A,B\}$ and $g=\{C,B\}$ in $\mathcal{F}^{(T|B)}$ with $|A\cap C|=a-1$ and $|B|=b$. It follows from Proposition\ \ref{pan2-1} (ii) and Remark\ \ref{pan2-2} that there exists a $b$-subset $D\subseteq[n]$ such that $f^\rho,g^\rho\in\mathcal{F}^{(T|D)}$. Additionally, if $f^\rho=\{G,D\}$ and $g^\rho=\{H,D\}$ then $|G\cap H|=a-1$. In a similar manner, $h^\rho\in\mathcal{F}^{(T|D)}$ in case when $h=\{X,B\}\in\mathcal{F}^{(T|B)}$ with $|X\cap A|=a-1$ or $|X\cap C|=a-1$. Along this idea of thought, Proposition\ \ref{pan2-6} implies that $h^\rho\in\mathcal{F}^{(T|D)}$ for every $h\in\mathcal{F}^{(T|B)}$. Therefore, $\Omega$ is a system of nontrivial blocks of $Aut(\Gamma(n,T))$ acting on $\mathcal{F}^{T}_n$.  \qed

\begin{lem}\label{pan3-5}\normalfont
Let $T=\{a,b\}\subseteq[n-1]$ such that $a<\frac{n}{2}<b<\frac{2n}{3}$ and $a+b+1\leq n$. Then $Aut(\Gamma(n,T))\cong S_n$.
\end{lem}
\demo By Proposition\ \ref{pan3-4}, it follows that $\Omega=\{\mathcal{F}^{(T|B)}_n: B\subseteq[n],|B|=b\}$ is a system of nontrivial blocks of $Aut(\Gamma(n,T))$ acting on $\mathcal{F}^{T}_n$. Consider the induced action $Aut(\Gamma(n,T))|_{\Omega}$ of $Aut(\Gamma(n,T))$ acting on $\Omega$. Let $\rho\in Aut(\Gamma(n,T))$ and ${\mathcal{F}}^{(T|A)}_n,{\mathcal{F}}^{(T|B)}_n,{\mathcal{F}}^{(T|C)}_n,{\mathcal{F}}^{(T|D)}_n\in\Omega$ such that ${{\mathcal{F}}^{(T|A)}_n}^{\rho}={{\mathcal{F}}^{(T|B)}_n}$ and ${{\mathcal{F}}^{(T|C)}_n}^{\rho}={{\mathcal{F}}^{(T|D)}_n}$. Note that $B\cup D=[n]$ if and only if $A\cup C=[n]$. Otherwise one of $(f,g),(f^\rho,g^\rho)$ is an edge and the other is not an edge will occur for some $f\in\mathcal{F}^{(T|A)}_n$ and $g\in\mathcal{F}^{(T|C)}_n$. Therefore, $Aut(\Gamma(n,T))|_{\Omega}\leq Aut(KG(n,n-b))$. On the other hand, the symmetric group on the set $[n]$ induces an automorphism group of $\Gamma(n,T)$, and thus $Aut(\Gamma(n,T))|_{\Omega}\cong Aut(KG(n,n-b))\cong S_n$. Thus, it suffices to check that no nonidentity automorphism can fix all the blocks of $\Omega$. Let $\rho$ be an automorphism in $Aut(\Gamma(n,T))$ such that ${{\mathcal{F}}^{(T|B)}_n}^{\rho}={{\mathcal{F}}^{(T|B)}_n}$ for all $B\subseteq[n]$ with $|B|=b$. Assume that there exist two distinct flags $f=\{A,B\}$ and $g=\{C,B\}$ in $\mathcal{F}^{(T|B)}_n$ such that $f^\rho=g$ for some $B\subseteq[n]$ with $|B|=b$. Since $f\neq g$, it follows that $A\neq C$, in other words, there exists an $i\in C\setminus A$. Hence, there exists a flag $h=\{D,E\}\in N(f)$ such that $i\in D$ and $|D|=a$ and $|E|=b$. By Fact\ \ref{pan2-0}, we see that $N(f)^\rho=N(g)$. However, it is clear that $N(g)\cap\mathcal{F}^{(T|E)}_n=\emptyset$, and which indicates that ${{\mathcal{F}}^{(T|E)}_n}^{\rho}\neq{{\mathcal{F}}^{(T|E)}_n}$, a contradiction.    \qed

\begin{lem}\label{pan3-6}\normalfont
Let $T=\{a,b\}\subseteq[n-1]$ such that $a\leq \frac{n}{3}-1$ and $b=\frac{2n}{3}$. Then the followings hold.

(i) $\Sigma=\{\mathcal{F}^{(T|A)}_n: A\subseteq[n],|A|=a\}$ and $\Omega=\{\mathcal{F}^{(T|B)}_n: B\subseteq[n],|B|=b\}$ are two systems of nontrivial blocks of $Aut(\Gamma(n,T))$ acting on $\mathcal{F}^{T}_n$.

(ii) $Aut(\Gamma(n,T))\cong S_n$.
\end{lem}
\demo Let $\rho$ be an automorphism in $Aut(\Gamma(n,T))$. Pick $f=\{A,B\}$ and $g=\{A,C\}$ in $\mathcal{F}^{(T|A)}_n$ such that $|B\cap C|=b-1$ and $|A|=a$. It follows from Proposition\ \ref{pan2-1} (iii) and Remark\ \ref{pan2-2} and Proposition\ \ref{pan2-4} and Proposition\ \ref{pan2-5} that there exists a $a$-subset $D\subseteq[n]$ such that $f^\rho,g^\rho\in\mathcal{F}^{(T|D)}$. An argument similar to the one used in the proof of Proposition\ \ref{pan3-2} shows that $\Sigma$ is a systems of nontrivial blocks of $Aut(\Gamma(n,T))$ acting on $\mathcal{F}^{T}_n$. In a similar manner, we can show that $\Omega$ is also a systems of nontrivial blocks of $Aut(\Gamma(n,T))$ acting on $\mathcal{F}^{T}_n$. Proceeding as in the proof of Lemma\ \ref{pan3-3} or Lemma\ \ref{pan3-5}, we have $Aut(\Gamma(n,T))\cong S_n$.   \qed

So far, we have not used \cite[Remark~5.13]{M} to give a positive answer to Question\ \ref{pan1-1}. On the other hand, by applying \cite[Lemma~2.1(b)]{M}, we obtain the following theorem.
\begin{thm}\label{pan3-7}\normalfont
Let $T=\{a,b\}\subseteq[n-1]$ with $a<\frac{n}{2}<b$ and $a+b\neq n$. Then $Aut(\Gamma(n,T))\cong S_n$.
\end{thm}

\section{Concluding Remarks}
Let $T=\{a,b\}\subseteq[n-1]$ with $a<b$. Note that three situations are left to consider, those are, $a+b=n$ and $b\leq\frac{n}{2}$ and $\frac{n}{2}\leq a$.

Let $T=\{a,b\}\subseteq[n-1]$ with $a<b$ and $a+b=n$. Suppose that $f=\{A,B\}$ and $g=\{C,D\}$ are two flags in $\mathcal{F}^{T}_n$ such that $|A|=|C|=a$ and $|B|=|D|=b$. Clearly, $(f,g)$ is an edge if and only if $C=\overline{B}$ and $D=\overline{A}$. Define $\Delta^{T(A|B)}_{n}=\{f,g\}$ where $f=\{A,B\}$ and $g=\{\overline{A},\overline{B}\}$ are in $\mathcal{F}^{T}_n$, and $\Omega=\{\Delta^{T(A|B)}_{n}:A\subset B\subset[n],|A|=a,|B|=b\}$. It is straightforward to see that $\Omega$ is a system of nontrivial blocks of $Aut(\Gamma(n,T))$ acting on $\mathcal{F}^{T}_n$, and further the following lemma holds.

\begin{lem}\label{pan4-1}\normalfont
Let $\Omega=\{\Delta^{T(A|B)}_{n}:A\subset B\subset[n],|A|=a,|B|=b\}$, where $T=\{a,b\}\subseteq[n-1]$ with $a<b$ and $a+b=n$. Then $Aut(\Gamma(n,T))\cong N\wr S_m$, where $N=\overbrace{S_2\times\cdot\cdot\cdot\times S_2}^{m}$ and $m=|\Omega|=\frac{1}{2}C^b_nC^a_b$.
\end{lem}

For the case that $b\leq\frac{n}{2}$, we consider a more general cases, as follows.

\begin{lem}\label{pan4-2}\normalfont
Let $\Sigma=\{\mathcal{F}^{(T|A)}_n: A\subseteq[n],|A|=t_r\}$, where $T=\{t_1,t_2,...,t_r\}\subseteq[n-1]$ such that $t_1<t_2<\cdot\cdot\cdot<t_r\leq\frac{n}{2}$. Then $\Sigma$ is a system of nontrivial blocks of $Aut(\Gamma(n,T))$ acting on $\mathcal{F}^{T}_n$. In addition, $Aut(\Gamma(n,T))\cong N\wr S_n$ where $N=\overbrace{S_m\times\cdot\cdot\cdot\times S_m}^{C^{t_r}_n}$ and $m=|\mathcal{F}^{(T|A)}_n|$ for a $\mathcal{F}^{(T|A)}_n\in\Sigma$.
\end{lem}

\demo Let $\rho\in Aut(\Gamma(n,T))$ and $\mathcal{F}^{(T|A)}_n\in\Sigma$. Suppose that $f\in\mathcal{F}^{(T|A)}$ and $f^\rho\in\mathcal{F}^{(T|B)}$ for some $\mathcal{F}^{(T|B)}\in\Sigma$. It suffices to prove that $g^\rho\in\mathcal{F}^{(T|B)}$ for any $g\in\mathcal{F}^{(T|A)}$. Proof by contradiction. Assume that there exists a $h\in\mathcal{F}^{(T|A)}$ such that $h^\rho\in\mathcal{F}^{(T|C)}\neq\mathcal{F}^{(T|B)}$. It is clear that $$N(f,h)=\bigcup\mathcal{F}^{(T|D)}_n,~{\rm{ where}}~ D\subseteq[n]\setminus A~ {\rm{with}}~ |D|=t_r.$$
Likewise, it is simple to see that
$$N(f^\rho,h^\rho)=\bigcup\mathcal{F}^{(T|D)},~{\rm{ where}}~D\subseteq[n]\setminus (B\cup C)~ {\rm{with}}~ |D|=t_r.$$
Obviously, $|N(f^\rho,h^\rho)|<|N(f,h)|$, and which is contradict to Fact\ \ref{pan2-0}. Therefore, $\Sigma$ is a system of nontrivial blocks of $Aut(\Gamma(n,T))$ acting on $\mathcal{F}^{T}_n$. 

Let $\rho\in Aut(\Gamma(n,T))$ and ${\mathcal{F}}^{(T|A)}_n,{\mathcal{F}}^{(T|B)}_n,{\mathcal{F}}^{(T|C)}_n,{\mathcal{F}}^{(T|D)}_n\in\Sigma$ such that ${{\mathcal{F}}^{(T|A)}_n}^{\rho}={{\mathcal{F}}^{(T|B)}_n}$ and ${{\mathcal{F}}^{(T|C)}_n}^{\rho}={{\mathcal{F}}^{(T|D)}_n}$. Clearly, $B\cap D=\emptyset$ if and only if $A\cap C=\emptyset$. Hence, we deduce that $Aut(\Gamma(n,T))|_{\Sigma}\leq Aut(KG(n,t_r))$. Additionally, the symmetric group on the set $[n]$ induces an automorphism group of $\Gamma(n,T)$, and therefore $Aut(\Gamma(n,T))|_{\Sigma}\cong Aut(KG(n,t_r))\cong S_n$. On the other hand, it is obvious that the kernel of the natural homomorphism from $Aut(\Gamma(n,T))$ to $Aut(\Gamma(n,T))|_\Sigma$ is the direct product of the symmetric groups on all blocks of $\Sigma$, and therefore $Aut(\Gamma(n,T))\cong N\wr S_n$ where $N=\overbrace{S_m\times\cdot\cdot\cdot\times S_m}^{C^{t_r}_n}$ and $m=|\mathcal{F}^{(T|A)}_n|$ for a $\mathcal{F}^{(T|A)}_n\in\Sigma$.   \qed

By \cite[Lemma~2.1(b)]{M} and Lemma\ \ref{pan4-2}, we derive the following corollary.

\begin{cor}\label{pan4-3}\normalfont
Let $T=\{t_1,t_2,...,t_r\}\subseteq[n-1]$ with $\frac{n}{2}\leq t_1<t_2<\cdot\cdot\cdot<t_r$. Then $Aut(\Gamma(n,T))\cong N\wr S_n$, where $N=\overbrace{S_m\times\cdot\cdot\cdot\times S_m}^{C^{t_1}_n}$ and $m=|\mathcal{F}^{(T|A)}_n|$ with $|A|=n-t_1$.
\end{cor}

\section{Acknowledgement}

We are very grateful to the anonymous referees for their useful suggestions and comments.

\end{document}